\documentclass[12pt]{article}
\usepackage{amssymb,amsmath, epsfig}
\include{epsf}

\usepackage{graphicx}
\usepackage{amsfonts}
\usepackage{amscd}
\usepackage[all]{xy}
\advance \textheight by \topskip
\makeatletter

\newcommand{\id}{{\mathrm{id}}}

\renewcommand{\theequation}{\thesection.\arabic{equation}}

\newcommand{\noi}{\noindent}
\newcommand{\cA}{{\cal {A}}}
\newcommand{\Gsub}{\underline{\G}}

\def\appendix#1
{
\addtocounter{section}{1}\setcounter{equation}{0}
 \renewcommand{\thesection}{\Alph{section}}
 \section*{Appendix~\thesection\protect\indent \parbox[t]{11.715cm}{#1}}
 \addcontentsline{toc}{section}{Appendix \thesection\ \ \ #1}
}

\newcommand{\be}{\begin{equation}}
\newcommand{\ee}{\end{equation}}
\newcommand{\beq}{\begin{equation}}
\newcommand{\eeq}{\end{equation}}
\newcommand{\bea}{\begin{eqnarray}}
\newcommand{\eea}{\end{eqnarray}}
\def\beqa{\begin{eqnarray}}
\def\eeqa{\end{eqnarray}}

\def\KK{{\mathbb K}}

\def\RR{{\mathbb R}}

\def\G{G}
\def\g{g}

\newcommand{\veps}{\varepsilon}

\newcommand{\cC}{{\cal {C}}}
\newcommand{\cB}{{\cal {B}}}

\newcommand{\eq}{\begin{equation}}
\newcommand{\eqa}{\begin{eqnarray}}
\newcommand{\en}{\end{equation}}
\newcommand{\ena}{\end{eqnarray}}

\newtheorem{theorem}{Theorem}[section]

\newtheorem{corollary}[theorem]{Corollary}

\newtheorem{definition}[theorem]{Definition}

\newtheorem{lemma}[theorem]{Lemma}

\newtheorem{proposition}[theorem]{Proposition}
\newtheorem{remark}[theorem]{Remark}

\def\KK{{\mathbb K}}

\def\RR{{\mathbb R}}

\newcommand{\prf}{{\noindent \bf Proof\; \; }}
\newcommand{\qed}{{\hfill $\Box$}}
\newcommand{\cH}{{\cal H}}

\begin{document}
\title{
\begin{flushright}
\end{flushright}
\vskip 1cm
{\bf Combinatorial Hopf algebraic description of the multiscale renormalization in quantum field theory } }
\author{
{\sf   Thomas Krajewski${}^{a}$\thanks{e-mail: thomas.krajewski@cpt.univ-mrs.fr},  
Vincent Rivasseau${}^{b}$\thanks{e-mail: vincent.rivasseau@th.u-psud.fr} 
 and
{\sf Adrian Tanasa}${}^{c}$\thanks{e-mail: adrian.tanasa@ens-lyon.org}}\\
}
\maketitle

\vskip-1.5cm

\vspace{2truecm}

\begin{abstract}
\noindent 
We define in this paper combinatorial Hopf algebras, on assigned Feynman graphs and on Gallavotti-Nicol\`o trees, 
 which we then prove to underly the multi-scale renormalization in quantum field theory.
Moreover,  morphisms between these Hopf algebras and the Connes-Kreimer Hopf algebras, 
on rooted trees and on 
Feynman graphs, are given. Finally, we show how this formalism can be used to investigate some algebraic properties of the effective expansion in multiscale renotmalization.
\end{abstract}




\tableofcontents

\section{Introduction and motivation}
\renewcommand{\theequation}{\thesection.\arabic{equation}}
\setcounter{equation}{0}
\label{introduction}

The interplay between {\it combinatorics} and {\it physics} has been recently very fruitful for both, and leads to the new and 
emerging interdisciplinary field of {\it combinatorial physics}.
For instance combinatorial tools have been successfully used for a better understanding of the algebraic structures underlying quantum mechanics (see \cite{flajolet}, \cite{5} and references within) and the interplay between combinatorics and statistical physics or integrable systems has been extensively studied both by combinatorists and by theoretical physicists.

In quantum field theory (QFT), a similar success is the elegant description of the combinatorial backbone of perturbative renormalization
{\it via} the combinatorial Connes-Kreimer Hopf algebra on Feynman graphs (see the original paper \cite{CK0} as well as section  $1.6$ of the book \cite{cm})\footnote{Let us also mention here that a Hopf algebraic description was also used to describe the combinatorics of perturbative renormalization on noncommutative Moyal space scalar QFT (where graphs are replaced by ribbon graphs, or combinatorial maps) 
\cite{io-fabien}, \cite{io-kreimer}; moreover, Connes-Kreimer-like Hopf algebras have been defined for quantum gravity spin-foam models \cite{fotini}, \cite{io-sf}.}.
The Connes-Kreimer Hopf algebra allows to recover the analytic expressions of a renormalized Feynman amplitude and the usual forest structure of the
subtraction operators, e.g. in the Bogoliubov-Parasiuk-Hepp-Zimmermann (BPHZ) renormalization, by using the recursive computation of the antipode,
which automatically generates all Zimmermann's forests with their correct weight.

This elegant point of view and its relationship with other mathematical problems such as the Riemann-Hilbert problem \cite{CK1,CK2} has made
renormalization a popular subject of mathematics. But there is a drawback: it has 
become so famous among mathematicians that it may have obscured for some of them the true physical meaning of renormalization.
Indeed the key physical notion in renormalization, namely the notion of scale, is absent or hidden in the Connes-Kreimer formalism. 

It is the goal of this paper to attract the attention of the mathematics community on this 
point and to propose a possible compromise by supplementing the Connes-Kreimer algebra with discrete scale assignments. The
corresponding algebra is generated by \emph{assigned graphs}, which are ordinary Feynman graphs supplemented with the assignment of
an integer to each edge. This integer physically represents the resolution scale of that edge or propagator. 

Indeed the modern version of renormalization, namely the renormalization group discovered by Wilson and followers \cite{Wil1},
tells us that the main purpose of renormalization is not to remove divergencies from Feynman amplitudes, nor to hide 
them into unobservable infinite bare parameters\footnote{This is particularly clear in the case of asymptoticallly free theories such as quantum chromodynamics, 
the theory of strong interactions, for which the bare coupling tends to \emph{zero}
at a high ultraviolet scale.}. Renormalization is much more general and powerful. It 
explains why and how, for physical systems with many coupled degrees of freedom, the \emph{laws of nature} \emph{change with the observation scale} \cite{Wil2}. 
This fundamental aspect of renormalization is captured mathematically by a \emph{multiscale analysis}. 

For general systems the Wilsonian slicing into scales can be implemented technically in many ways (block-spins, wavelet analysis, etc...).
The most convenient technique in the context of QFT \emph{slices} the propagator of the theory according to
a geometric sequence of cutoffs. 
Each slice represents a particular energy scale, and has a particular
spatial resolution power; it has an ultraviolet and an infrared cutoff with constant ratio between both. 
The renormalization group then performs many times the same step, namely functional integral over a slice or fluctuation field and 
computation of the resulting effective action for the remaining sum of fields of lower slices, called the background field \cite{Wil2}.

The need for such a discrete multiscale analysis of QFT
was quite independently also discovered by mathematical physicists such as J. Glimm, A. Jaffe and their followers
of the constructive field theory program \cite{Erice}. They called it the \emph{phase space expansion}.
Over the years this constructive program, in which perturbative QFT is summed, 
effectively  merged completely with the Wilsonian renormalization group 
approach and became its mathematically rigorous version\footnote{Let us remark at this point that 
(contrary to a belief sometimes heard in the mathematics community)
the residues or the individual renormalized Feynman amplitudes do not correspond to any physical observables in
QFT. Indeed any measure always involve not a single Feynman amplitude but the infinite sum of
such amplitudes compatible with a given set of external legs hence also with a certain \emph{resolution power}. Only the value of this infinite sum, which is what  constructive theory is after, has physical meaning.}. Constructive analysis comes at a high price: many elegant perturbative tools in QFT
such as dimensional regularization, dimensional renormalization, and differential (rather than finite-difference) renormalization group equations, had to be discarded
by the constructive community in favor of \emph{discrete} multiscale analysis, which remains up to now the only tool with 
proven constructive power. For a general presentation of these 
views and of multiscale renormalization, 
see \cite{book-rivasseau}.
 
Returning 
to the more limited and specific context of perturbative renormalization of Feynman amplitudes,
multiscale analysis was first developed systematically in \cite{FMRS} and
\cite{GN}. Initially these authors were motivated by the desire to understand and simplify the proof of uniform bounds on renormalized amplitudes 
implying ``local Borel summability" \cite{dCR}, which had been soon followed also by the construction and Borel summability of \emph{planar} 
asymptotically free renormalized theories, such as "wrong sign" planar $\phi^4_4$ \cite{tHooft, Riv}.

Multiscale analysis evolved over the years into a versatile technique to understand and 
analyze renormalization and the renormalization group in new contexts.
It was suitably generalized to the condensed matter case in which the Fermionic 
propagator is sliced in a sequence
of scales pinching closer and closer the Fermi surface \cite{BG1,FT1}. 
This technique provided the backbone for the rigorous analysis
of correlated quantum Fermions at low temperature, such as Fermi and Luttinger liquids in 
one, two and three spatial dimensions (see \cite{Rivasseau:2011ri} and the many references therein for a recent review of 
this large, active and mature field of mathematical physics). 

More recently the multislice analysis has been used to prove perturbative renormalizability at all orders for
radically new quantum field theories in which the interaction is non-local and the usual intuition of zero momentum subtraction
around local parts fails. Such new models include the first examples of renormalizable
noncommutative quantum field theories (see \cite{GW}, \cite{propa}, \cite{4men}, \cite{gn} and \cite{GMRT})
and of tensor group field theories \cite{BGR,COR}. 
The latter models might be relevant for the long term goal of 
quantization of gravity \cite{Rivasseau:2011hm}, but also for the more concrete 
analysis of statistical physics in random geometry \cite{Bonzom:2011ev} or with long range interactions, such as spin glasses \cite{BGS}.

Multiscale analysis is characterized by the fact that 
the contraction and subtraction operations that implement renormalization are not effectuated blindly. 
They make physical sense only for so-called \emph{high}  subgraphs, i.e. connected subgraphs 
which have all their internal scales higher than any of their external scales. It is solely for such subgraphs that the comparison of their amplitudes to a local part 
makes sense\footnote{Physically this is nothing but the trivial observation that objects with a certain size look local only when observed through probes that do not distinguish their internal structure.}. It is this distinction which in turns launches the renormalization group flow, hence the motion of effective constants
with scale. Assigned graphs allow to define such high parts, whether general graphs do not; hence we feel they should become part of the combinatoric Hopf algebra framework used by mathematicians to describe renormalization.

\medskip 

In this paper, we therefore define a new Hopf algebra which is meant to describe the combinatorial soul of this discrete multiscale renormalization technique.
In order to do that we define  {\it assigned graphs}  as Feynman graphs together with a scale assignment of their edges. 
The desired combinatorial Hopf algebra is then defined on the space freely generated by these assigned graphs.  The coproduct has then to take into account the supplementary scale information of the assigned graphs: one only sums over the particular class of {\it high subgraphs}. For example, in the case of 
the scalar $\phi^4$ model, one does not need to sum over {\it all} subgraphs with two- or four-external edges, as is done in the standard definition of the Connes-Kreimer 
coproduct.

\medskip

Let us also mention that 
in this paper we deal with the $\phi^4$ model, even though our results can generalize in a straightforward manner to 
more general renormalizable QFTs.

\section{Feynman graph expansion and multiscale renormalization}

In this section, we give a short overview of Feynman graphs and multiscale expansions in quantum field theory.

\subsection{From path integral to Feynman graphs}

In its most general acceptance, QFT can be defined as the study of quantum (or stochastic) dynamical systems involving continuous degrees of freedom. In the euclidian path integral approach, one has to define the path integral representing the expectation value of an observable, heuristically written as
\begin{equation}
\langle{\cal O}\rangle=\int\frac{d\mu(\phi)\,{\cal O}[\phi]\,\exp-{\cal S}[\phi]}
{\int d\mu(\phi)\,\exp-{\cal S}[\phi]}
.\end{equation}
The integration is over a suitable space of fields $\phi:\,{\Bbb R}^{D}\rightarrow{\Bbb R}$, ${\cal S}[\phi]$ is the action and ${\cal O}$ an observable. 

In the simplest case (the so called $\phi^{N}_{D}$ euclidian field theory), the action can be written as
\begin{equation}
{\cal S}[\phi]=\int_{{\Bbb R}^{D}} d^{D}x\, \Big\{ \frac{1}{2}\phi(x)\big(-\Delta+m^{2}\big)\phi(x)+\frac{\lambda}{N!}\phi(x)^{N}\Big\}
\end{equation}
with $\Delta$ the Laplacian and $m$ and $\lambda$ two positive real numbers, identified with the mass and the coupling constant of the theory. The observable are usually taken to be products of the fields at different space-time points, ${\cal O}[\phi]=\phi(x_{1})\dots\phi(x_{n})$, whose expectation value define the $n$-point correlation functions.

In the free field case $\lambda=0$, the path integral is Gau\ss ian and is readily computed using Wick's theorem. With a suitably normalized measure, the expectation value of a product of fields reads
\begin{equation}
\int d\mu_{C}(\phi) \phi(x_{1})\dots\phi(x_{n})=
\left\{
\begin{array}{cl}
0&\text{if }n \text{ is odd}\\
{\displaystyle \sum_{\text{pairings of}\atop \left\{ 1,2,\dots,n\right\}}}C(x_{i_{1}},x_{i_{2}})\cdots C(x_{i_{\frac{n}{2}\!-\!1}},x_{I_{\frac{n}{2}}})&\text{if }n \text{ is even}
\end{array} 
\right.
\end{equation}
where the covariance $C(x,y)$ is the kernel of the inverse of $-\Delta+m^{2}$.

In the non Gaussian case, we expand the integrand as a formal power series in $\lambda$ and perform all the integrals using Wick's theorem. 
Each term we obtain this way is called a Wick contraction. Wick contractions naturally define graphs, called Feynman graphs.  Collecting Wick contractions that correspond to the same graph, we obtain an expansion of the correlation functions as a sum over graphs,
\begin{equation}
\int d\mu_{C}\phi(x_{1})\dots\phi(x_{n})=\sum_{G\text{ graph with}\atop n\text{ external edges}}
\frac{\lambda^{v(G)}}{\sigma(G)}A(G)[x_{1},\dots,x_{n}]
\label{pathintegral}
\end{equation}
The Feynman graphs have $n$ labelled univalent vertices (associated to the variables $x_{1},\dots,x_{n}$ and $v(G)$ $N$-valent vertices corresponding to the interaction monomial $\phi^{N}(x)$. The edges related to two $D$-valent vertices are called internal edges and the other edges are called external. Because of the variables  $x_{1},\dots,x_{n}$, there are labels on the external edges, while the internal ones are unlabeled so that we are summing over isomorphism classes of graphs with fixes external edges . This is accounted for by the symmetry factor $\sigma(G)$ defined as follows. 

When expanding the path integral, we label the vertices and for each vertex, we also label the half-edges emanating from it, so that the half-edges are labelled by pairs $(v,p)$ . Then, a Wick contraction is just a partition of the indices of half lines into pairs and the symmetry factor $\sigma(G)$ is the subgroup of the group of permutations of all these labels of the internal lines that preserve this partition, with $G$ the corresponding isomorphism class.  In the general case, there are $v(G)! (N!)^{v(G)}$ labelings of the internal half lines,  so that there are $\frac{v(G)! (N!)^{v(G)}}{\sigma(G)}$ Wick pairings associated to a given isomorphism class.

From an analytic point of view, one has to remember that the kernel $C(x,y)$ is a distribution and the Feynman graph amplitudes are not well defined since they involve products of distributions. In order to over come this problem, one first regulates the theory, replacing the distribution by some function $C_{\rho}(x,y)$ depending on a regulator $\rho$, in such a way that we recover $C(x,y)=\lim_{\rho\rightarrow\infty} C_{\rho}(x,y)$. Then, a recursive operation is performed on the Feynman graphs amplitudes in such a way that they are well defined in the limit $\rho\rightarrow\infty$, for the so-called  renormalizable theories $\phi_{4}^{4},\phi_{6}^{3},\phi_{3}^{6},\dots$. For graphs without subdivergent graphs, this operation is additive but otherwise one has to first renormalize the subdivergencies. The corresponding operation is polynomial and involves a sum over all the forests of the graphs, as will be made more precise later. We refer to \cite{Collins} for a detailed overview of perturbative 
renormalization and the Boliubov-Parasiuk-Hepp-Zimmermann (BPHZ) forest formula. 

Even if very successful, the BPHZ forest formula has an important drawback: it does not implement Wilson's idea that path integrals must be computed by first integrating over small distance degrees of freedom. To implement this idea, it is convenient to use multiscale analysis.

\subsection{Multi-scale renormalization in a nutshell}
\renewcommand{\theequation}{\thesection.\arabic{equation}}
\setcounter{equation}{0}
\label{sec:ms}

As mentioned in the introduction, the multiscale analysis of renormalization which is at the core of the Wilsonian approach
to relies on a geometrically growing sequence of discrete scales.
There are two main technical ways to create the sequence of  scales:
\begin{itemize}
\item 
block spinning of the field variables in direct space, that is defining $\phi = \phi_f + \phi_b$, where 
the background field $\phi_b $ is the local \emph{average} of $\phi$
with respect to a lattice of cubes of side size $M$, and $\phi_f$, the fluctuation field is simply the difference
between the field and the background field;

\item
slicing the propagator $C$ as $C_f + C_b$, where 
$C_f$ has both infrared and ultraviolet cutoff with fixed ratio $M$
and $C_b$ has only an ultraviolet cutoff, which is the infrared cutoff of $C_f$; 
in that case the slicing induces an orthogonal decomposition
of the field as $\phi = \phi_f + \phi_b$, where $\phi_f$ is distributed according to $C_f$
and $\phi_b$ according to $C_b$.
\end{itemize}

The first technique is more general and can apply to any statistical mechanics system 
but requires a discretization through lattices.
The second technique is the most elegant and clearly best adapted to perturbative renormalization
theory around a propagator with non-trivial spectrum. More precisely an excellent compromise for propagator with a 
positive spectrum is the parametric slicing:

\begin{definition}[Parametric Slicing]
Let $C= 1/H$ be the propagator of the theory. The parametric slicing is
\bea  C &=& \int_0^\infty  e^{- \alpha H}  d \alpha \; ,
\;   \sum_{i=0}^{\infty}  C^{i}   \label{decoab1}  \\
C^{i} &=&  \int_{M^{-2i}}^{M^{-2(i-1)}} e^{- \alpha H}  d \alpha \; ,
\; C^{0} =  \int_{1}^{\infty} e^{- \alpha H}  d \alpha .
\eea
\end{definition}
The natural ultraviolet cutoff on the theory is then
\beq  C_{\rho} =  \sum_{i=0}^{\rho}  C^{i}   \label{decocut1} 
\eeq
for finite and large integer $\rho$. In the case of the Laplacian plus mass on $\RR^d$ we get the following slices
\bea   \label{deco1}
C^{i} &=&  \int_{M^{-2i}}^{M^{-2(i-1)}} e^{-m^{2}\alpha -  
{\vert x-y \vert^{2} \over 4\alpha }}  {d\alpha \over \alpha^{d/2}}   \label{deco2} \\
 C^{0} &=&  \int_{1}^{\infty} e^{-m^{2}\alpha -  
{\vert x-y \vert^{2} \over 4\alpha }}  {d\alpha \over \alpha^{d/2}} .  \label{deco3}
\eea
$\alpha$ being dual to $p^{2}$, one should consider each propagator  $C^{i}$
as corresponding to a theory with
both an ultraviolet and an infrared cutoff. They differ by the fixed 
multiplicative constant $M$, the momentum slice ``thickness". 

The decomposition \eqref{deco1}-\eqref{deco3} is the multislice representation. 
From the general definition of Gaussian measures follows an
associated decomposition of the Gaussian measure $d\mu_{\rho}$ of covariance $C_{\rho}$ 
into a product of independent Gaussian measures 
$d\mu^{i}$ with covariance $C^i$. Similarly the random field $\phi_{\rho}$ distributed 
according to $d\mu_{\rho}$ is the sum of \emph{independent} random variables 
$\phi^{i}$ distributed according to $d\mu^{i}$:

\beq \phi_{\rho} = \sum_{i=0}^{\rho} \phi^i ;\quad d\mu_{\rho}(\phi_{\rho})
= \otimes_{i=0}^{\rho} d\mu^{i}(\phi^{i})    \eeq

This independentness of the fields at each scale in turns leads in the perturbative analysis 
of the corresponding functional integral to a sum over \emph{assigned graphs}, that is graphs which have an
integer associated to each edge, namely its scale.

\begin{definition}
 A {\bf scale assignment} $\mu$ for a Feynman graph with labelled internal edges is a list of positive integers 
 $i_\ell$, $\ell=1,\ldots, E$ associated to the internal edges of the respective Feynman graph (where $E$ is the number of internal edges of the graph).
\end{definition}

\begin{figure}[h]
\centerline{   \includegraphics[width=7cm]{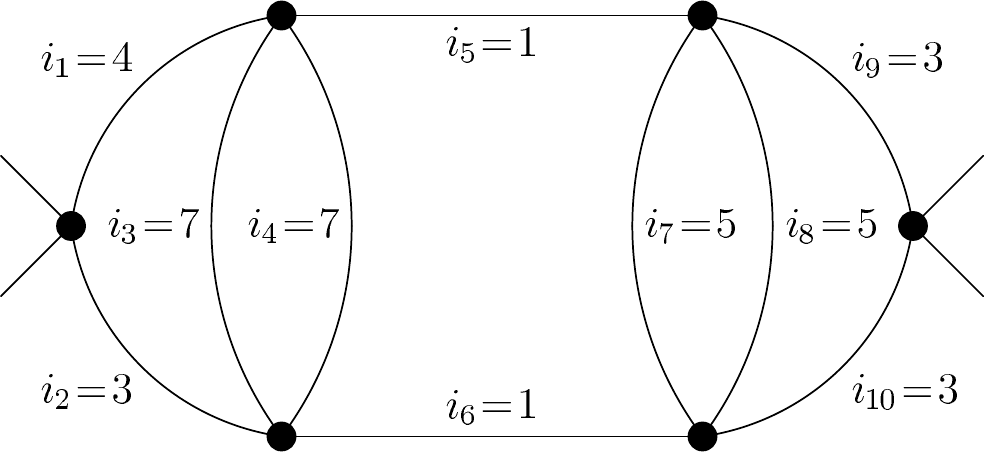}   }
\caption{A Feynman graph with a scale assignment; it has $10$ internal edges and $4$ external edges.}
\label{scales1}
\end{figure}

Let us emphasize here that the integers in the Definition above are bounded by the discrete cutoff $\rho$. One further has:

\begin{definition}
An {\bf assigned graph} $(G, \mu)$ is an isomorphism class of couples formed by the one particle irreducible (1PI)  edge labeled Feynman graph $G$, together with a scale assignment $\mu$.
\end{definition}

In physics, a 1PI graph is a graph that cannot be disconnected by cutting an arbitrary line. It is also called 2-edge connected in the mathematical literature.

\begin{remark}
Assigned graphs correspond to graphs whose edges are labeled by the scales. They can be seen as a particular class of decorated graphs.


\end{remark}

\begin{definition}
An  {\bf assigned subgraph} $(g,\nu)$ of a given assigned graph $(G,\mu)$ is constructed in the following way (see previous section). One considers a subgraph $g$ of the Feynman graph $G$, in the usual QFT way. The scale assignment $\nu$ of $G$ is given by the restriction of the scale assignment $\mu$ to the internal edges of $g$ (which are also internal edges of $G$). Moreover, the external edges of $g$ which are internal edges of $G$ have the scale assignment which is attributed to them by $\mu$. The same holds for the external edges of $g$ which are external edges of $G$.
\end{definition}

Furthermore, one can define the usual graph theoretical notions (number of edges, vertices, (independent) loops {\it etc.}) of an assigned graph $(G,\mu)$ as the respective notions of the Feynman graph $G$. Moreover, we call $(G,\mu)$ an $N-$point assigned graph if $G$ is an point $1$PI Feynman graph.

We then define the internal and external index for a subgraph $(g, \nu)$ of an assigned graph $(G, \mu)$ as:
\beq   i_g(\mu)  = \inf\limits_{l \in g}  \mu(l)  \label{II.1.10}
\eeq
\beq   e_g(\mu)  = \sup\limits_{l {\rm \ external \ line \ of \ } g}  
\mu(l)  \label{II.1.11}
\eeq
(with the $\mu$ dependence sometimes omitted for shortness). 

\begin{definition}
Let $(G,\mu)$ an assigned graph. We say that a subgraph $(g,\nu)$ is a {\bf high subgraph} if
\begin{itemize}
\item $g$ is connected

\item the internal index of $g$ is higher than its external index:
\beq   e_g(\mu)  < i_g(\mu)  \quad  \quad {\rm (high \ condition)}.  
\label{II.1.12}
\eeq

\end{itemize}
\end{definition}

One can  associate to a  connected  assigned graph $(G,\mu)$,
the Gallavotti-Nicol\`o tree $T_{(G,\mu)}$, which is defined in the following way (see the book  \cite{book-rivasseau} for more details on this):. 

\begin{definition}
\label{def-GN}
The Gallavotti-Nicol\`o tree $T_{(G,\mu)}$  is a rooted tree whose nodes at a distance $i$ from the root are decorated with the connected high subgraphs $G_{c}^{i}$ with scales $\geq i$ and whose arrows join the nodes decorated with $G^{i}_{c}$ and $G_{c'}^{i-1}$ if and only if $G_{c}^{i}$ is a high subgraph of $G_{c'}^{i-1}$.
\end{definition}

\begin{figure}[h]
\begin{tabular}{cc}
 \includegraphics[width=6.5cm]{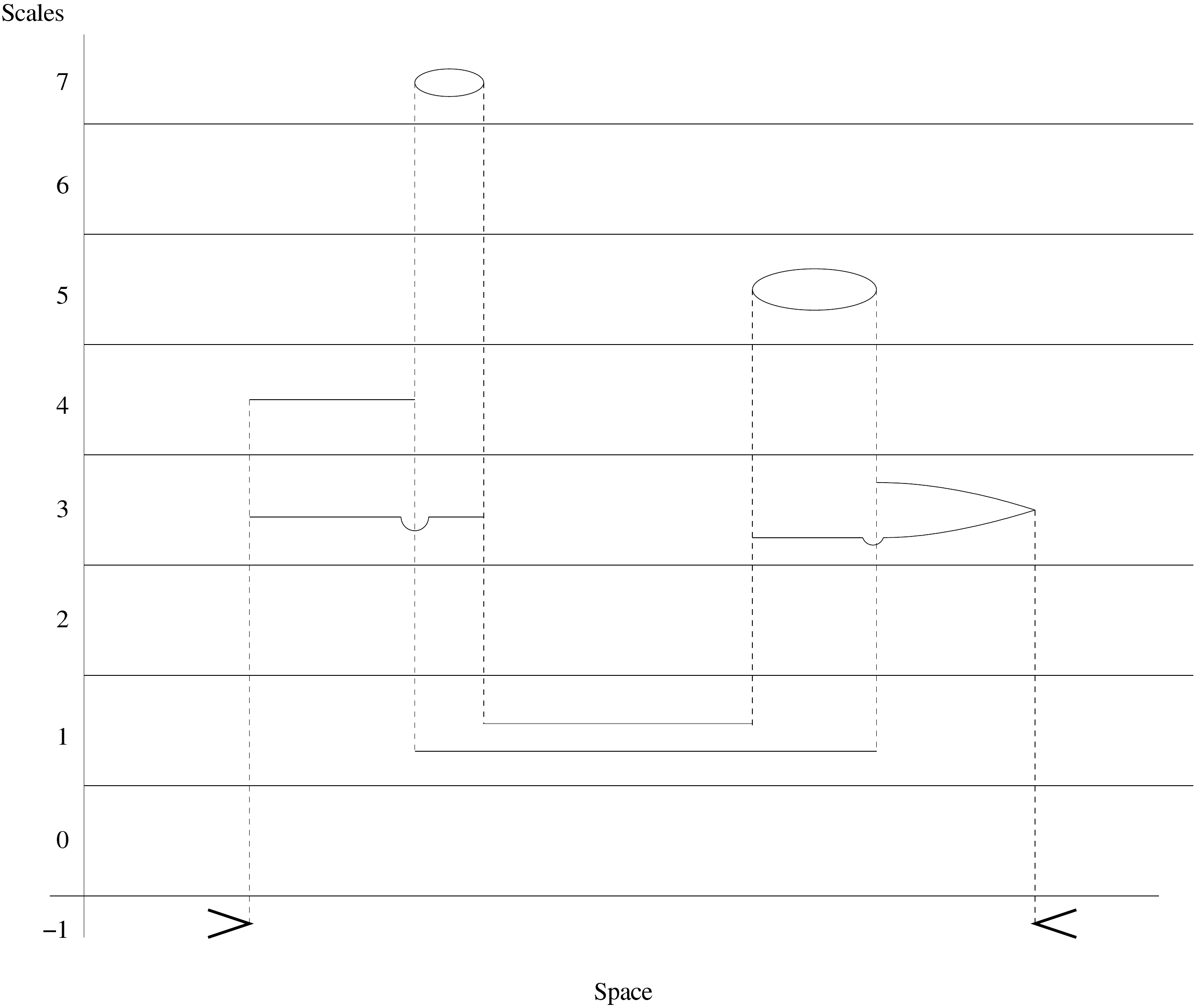} & \includegraphics[width=6.5cm]{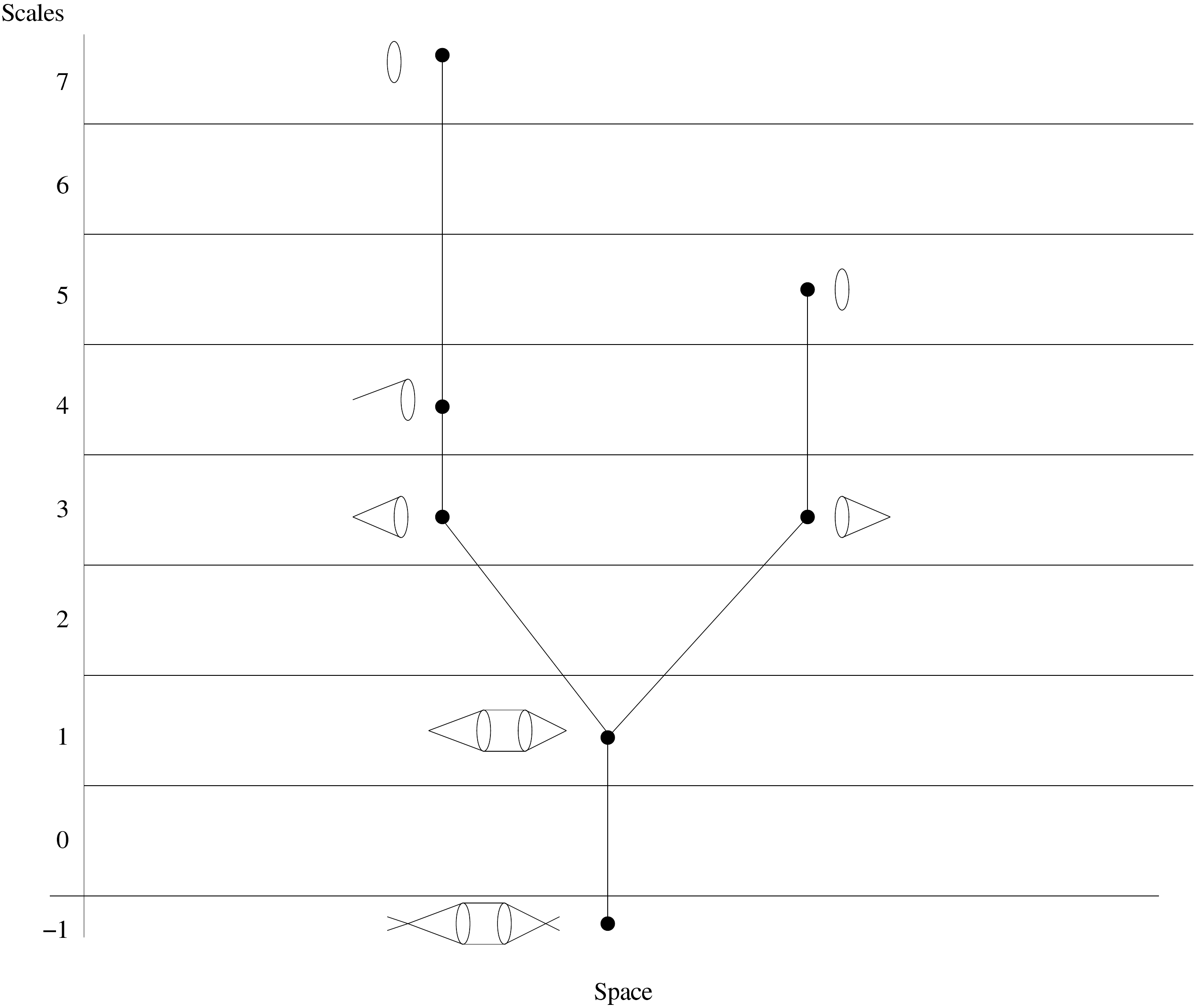} \\
 {\small a) Phase space representation.}&{\small b) Gallavotti-Nicol\`o tree.}
\end{tabular}

\caption{Phase space representation of the graph of Fig. \ref{scales1} and the associated 
Gallavotti-Nicol\`o tree.}
\label{scales3}
\end{figure}

In order to represent the Gallavotti-Nicol\`o tree, it is convenient to adopt a phase-space representation with positions on the horizontal axis and scales on the vertical one.  The phase space representation of the graph of figure \ref{scales1} and its 
Gallavotti-Nicol\`o tree can be found on figure  \ref{scales3}.

High subgraphs are partially ordered by inclusion. An essential result is that they form a
(Zimmermann's) forest in the following sense 
\begin{lemma}
Let $(G,\mu)$ be a fixed graph and scale assignment.
The set of high subgraphs is a forest, in the sense that if $g_1$ and $g_2$ are both high, we have either $g_1 \subset g_2$,
or $g_2 \subset g_1$ or $g_1 \cap g_2 = \emptyset$ . 
\end{lemma}
\prf
Suppose we could find $S_1$ and $S_2$
with a non trivial intersection; in this case $S_1$ would have an external line
$S_2$ and conversely; but the scale of any of these two lines should be both strictly larger
and strictly smaller than the other, which is impossible.
\qed

Usually the final graph $G$ is connected and this inclusion forest of  high subgraphs forms a tree which is nothing but the celebrated 
"Gallavotti-Nicol\`o" tree \cite{GN}.

\section{Hopf algebras on assigned graphs and the combinatorics of multi-scale renormalization}
\renewcommand{\theequation}{\thesection.\arabic{equation}}
\setcounter{equation}{0}
\label{sec:hopf}

{\subsection{Some algebra}}
\label{sec:alg}
\renewcommand{\theequation}{\thesection.\arabic{equation}}
\setcounter{equation}{0}


In this section we briefly recall, following \cite{io-fabien},  the definitions of the algebraic notions that will be used in the sequel.

\begin{definition}[Algebra]\label{defn:algebra}
  A unital associative algebra $\cal A$ over a field $\mathbb K$ is a $\mathbb K$-linear space endowed with two algebra homomorphisms: 
  \begin{itemize}
  \item a product $m :\cA\otimes\cA\to\cA$ satisfying the \emph{associativity} condition:
    \begin{align}
      m\circ(m\otimes\id)(\G)=&m\circ(\id\otimes m)(\G),
\ \forall \G\in \cA^{\otimes\, 3},
\label{eq:asso}
    \end{align}
  \item a unit $u :\KK\to\cA$ satisfying:
    \begin{align}
      m\circ(u\otimes\id)(1\otimes\G)=&\G=m\circ(\id\otimes u)(\G\otimes 1),\ \forall \G\in \cA.
    \end{align}
  \end{itemize}
\end{definition}



\begin{definition}\label{defn:coalgebra}
  A (coassociative, counital) {\bf coalgebra} $\cC$ over a field $\KK$ is a $\KK$-linear space endowed with two linear homomorphisms: 
  \begin{itemize}
  \item a coproduct $\Delta :\cC\to\cC\otimes\cC$ satisfying the \emph{coassociativity} condition:
    \begin{align}
      \forall\G\in\cC,\,(\Delta\otimes\id)\circ\Delta(\G)=&(\id\otimes \Delta)\circ\Delta(\G),\label{eq:coasso}
    \end{align}
  \item a counit $\veps :\cC\to\KK$ satisfying:
    \begin{align}
      \forall \G\in\cC,\,(\veps\otimes\id)\circ\Delta(\G)=&\G=(\id\otimes\veps)\circ\Delta(\G).    \end{align}
  \end{itemize}
\end{definition}

\begin{definition}
  A {\bf bialgebra} $\cB$ over a field $\KK$ is a $\KK$-linear space endowed with both an algebra and a coalgebra structure (see Definitions \ref{defn:algebra} and \ref{defn:coalgebra}) such that the coproduct and the counit are unital algebra homomorphisms (or equivalently the product and unit are coalgebra homomorphisms):
  \begin{subequations}
    \label{eq:compatibility}
    \begin{align}
      \Delta\circ m_{\cB}=&m_{\cB\otimes\cB}\circ(\Delta\otimes\Delta),\ \Delta (1_\cB)=1_\cB\otimes 1_\cB,\\
      \veps\circ m_{\cB}=&m_{\KK}\circ(\veps\otimes\veps),\ \veps (1_\cB)=1.
      \end{align}
    \end{subequations}
\end{definition}

\begin{definition}
  A {\bf graded bialgebra} is a bialgebra graded as a linear space: 
  \begin{align}
      \cB=\bigoplus_{n=0}^\infty\cB^{(n)}   
  \end{align}
  such that the grading is compatible with the algebra and coalgebra structures:
  \begin{align}
    \cB^{(n)}\cB^{(m)} \subseteq \cB^{(n+m)}\text{ and }\Delta\cB^{(n)}\subseteq\bigoplus_{k=0}^n\cB^{(k)}\otimes\cB^{(n-k)} ..
  \end{align}
\end{definition}

\begin{definition}
  A {\bf connected bialgebra} is a graded bialgebra $\cB$ for which $\cB^{(0)}=u(\KK)$.
\end{definition}


\begin{definition}
  A {\bf Hopf algebra} $\cH$ over a field $\KK$ is a bialgebra over $\KK$ equipped with an antipode map $S:\cH\to\cH$ obeying:
  \begin{align}
    m\circ(S\otimes\id)\circ\Delta=&u\circ\veps=m\circ(\id\otimes S)\circ\Delta.
  \end{align}
\end{definition}

We now end this section by recalling the following result:

\begin{lemma}[\cite{manchon}]
\label{lema-manchon}
  Any connected graded bialgebra is a Hopf algebra whose antipode is given by $S(1_\cH)=1_\cH$ and recursively by any of the two following formulas for $\G\neq 1_\cH$:
  \begin{subequations}
    \begin{align}
      S(\G)=&-\G-\sum_{(\G)}S(\G')\G'',\label{eq:Srecurs}\\
      S(\G)=&-\G-\sum_{(\G)}\G'S(\G'')
    \end{align}
  \end{subequations}
where we used Sweedler's notation.
\end{lemma}
It turns out that commutative Hopf algebras naturally give to a group structure on the space of characters.

\begin{definition}
A character of a commutative Hopf algebra is a linear map $\alpha$ from ${\cal H}$ to the ground field ${\Bbb K}$ such that $\alpha(\G_{1}\G_{2})=\alpha(\G_{1})\alpha(\G_{2})$ for any $\G_{1},\G_{2}\in{\cal H}$.  
\end{definition}

The group structure on the set of characters is given by the convolution product.

\begin{proposition}
The set $G$ of characters of ${\cal H}$ is a group for the multiplication law
\begin{equation}
\alpha\ast\beta=(\alpha\otimes\beta)\Delta
\end{equation}
with inverse $\alpha^{-1\ast}=\alpha\circ S$ and unit $\epsilon$.
\end{proposition}

For graded connected  Hopf algebras, characters form a Lie group whose Lie algebra are made of infinitesimal characters, defined as follows.

\begin{definition}
An infinitesimal character $\delta$ is a linear map from ${\cal H}$ to ${\Bbb K}$ such that $\delta(\G_{1}\G_{2})=\epsilon(\G_{1})\delta(\G_{2})+\delta(\G_{1})\epsilon(\G_{2})$ for any $\G_{1},\G_{2}\in{\cal H}$
\end{definition}

Infinitesimal characters define a Lie algebra ${\cal G}$ which is the Lie algebra of $G$. The convolution exponential
\begin{equation}
\alpha=\exp_{\ast}(\delta)=\sum_{n}\frac{\overbrace{\delta\ast\cdots\ast\delta}^{n\,\text{times}}}{n!}
\quad\Leftrightarrow\quad
\delta=\log_{\ast}=\sum_{n\leq 1}(-1)^{n-1}\frac{\overbrace{(\alpha-\epsilon)\ast\cdots\ast(\alpha-\epsilon)}^{n\,\text{times}}}{n}
\end{equation}

The Hopf algebra ${\cal H}$ can be understood as the algebra of functions from $G$ to ${\Bbb K}$.

For further details on this topic, the interested reader  can refer
for example to \cite{Kassel} or \cite{Dascalescu}.

\subsection{Hopf algebra structures on the Gallavotti-Nicol\`o trees}
\renewcommand{\theequation}{\thesection.\arabic{equation}}
\setcounter{equation}{0}
\label{sec:hopft}

One has:

\begin{proposition}
Let $T_{(G,\mu)}$ be the Gallavotti-Nicol\`o tree associated with the assigned graph $(G,\mu)$.
\begin{enumerate}
\item The root of $T_{(G,\mu)}$ is decorated with $G$ itself.
\item The leaves of  $T_{(G,\mu)}$ all are at distance $\rho$ from the root, with $\rho$ the ultraviolet cutoff.
\item
If the scales $j$ and $k$ are such that the scales $i$ obeying $j\leq i\leq k$ do not appear in $(G,\mu)$, then any subtree of $T_{(G,\mu)}$ whose root is a at distance $j$ and leaves at a distance $k$  from the root of $T_{(G,\mu)}$ does not branch and has all its nodes decorated by the same graph.
 \end{enumerate}
\end{proposition}
{\it Proof:} The items above follow as a direct consequence of the Definition \ref{def-GN} of the Gallavotti-Nicol\`o trees. (QED)

\medskip

In order to define the Hopf algebra underlying multiscale renormalization 
on the Gallavotti-Nicol\`o trees, it is useful to introduce the following terminology. If $T' \subset T_{(G,\mu)}$ is a subtree, we define its completion $\overline{T'}=T_{(G',\mu')}$ to be the 
Gallavotti-Nicol\`o tree  associated with its root $(G',\mu')$. Furthermore, we define an admissible cut $C$ to be a non empty subset of $|C|$ arrows of $T_{(G,\mu)}$ that join nodes decorated by two different graphs, the graph farther form the root having two of four external edges and such that any path form the leaves to the root contains at most one arrow in $C$. Removing the arrows in $C$ we get a subtree $T_{<}$ that contain the root and trees $T_{>}^{n}$ that do not contain the root.

\begin{proposition}
The free commutative algebra ${\cal H}_{\mathrm{GN}}$ generated by all Gallavotti-Nicol\`o trees is a graded Hopf algebra whose counit and coproduct are defined on the generators by  $\epsilon(T)=0$ and
\begin{equation}
\Delta(T)=T\otimes 1+1\otimes T+\sum_{C\atop\mathrm{admissible \,cut}}\Big(\prod_{1\leq n\leq |C|}\overline{T^{n}_{>}}\Big)\otimes \overline{T_{<}}.
\end{equation}
Its grading  is $n(T)=\#\{\mbox{arrows joining nodes decorated with different graphs}\}+1$.
\end{proposition}
{\it Proof.}  
The only non trivial assertions to check are the coassociativity of the coproduct and the existence of the grading and the antipode. The first proof is analogous to the proof of the coassociativity of the coproduct in the algebra of rooted trees, see \cite{CK0}. The assertion pertaining to the grading is easy to check as any cut reduces the number of of arrows joining nodes with different graphs by the number of cut edges. Finally, for any graded commutative bigebra there is a recursive construction of the antipode, as given in \cite{fgb}. (QED)


\subsection{Operations on assigned graphs}
\renewcommand{\theequation}{\thesection.\arabic{equation}}
\setcounter{equation}{0}
\label{sec:ag}

In this section we define several operations which we need in the rest of the paper.

\medskip

We now define $\underline{(G,\mu)}$ to be the set of assigned graphs formed by high subgraphs of  the assigned graph ${(G,\mu)}$ whose connected components are 1PI have two or four edges external edges. An external edge of a subgraph is an edge of $G$ attached to a vertex in $g$ which is not an internal edge of $g$.

\begin{definition}
\label{def:shrink}
The {shrinking} of a two- or a four-point assigned subgraph $(g,\nu)$ inside an assigned graphs $(G,\mu)$ is defined in the following way. The shrinking of the subgraph $g$ inside the Feynman graph $G$ is done in the usual QFT way, {\it i. e.} the subgraph is replaced by a vertex 
(the internal structure of $g$ vanishes); one has the cograph $G/g$.
The scale assignment $\mu/\nu$ of the cograph $G/g$ is given by the initial scale assignment
$\mu$, where we have erased the scale assignment of the internal edges of $g$ (if two external edges are added when shrinking a two-point function, they are assigned a non-dangerous integer). We call the resulting assigned graph $(G/g,\mu/\nu)$ a {assigned cograph}.
\end{definition}

\begin{remark}
The shrinking operation corresponds to the  wave function or mass renormalization, for a two-point subgraph, or to the coupling constant renormalization for a four-point subgraph. In the case of a wave-function renormalization, a decoration indicating the two derivative couplings of the Laplacian 
must be added to the shrinked two point vertex  to distinguish that renormalization from the mass renormalization.
\end{remark}

\begin{definition}
\label{def:gluing}
The {gluing data} $\circ$ for the insertion of a two- and respectively four-point assigned graph $(g,\nu)$ into the propagator and respectively the vertex of an assigned graph $(G,\mu)$ is 
given by a bijection between the external edges of $g$ and the two half-edges of the propagators or respectively the four half-edges of the vertex. It is
defined only if the external assignment indices for $(g,\nu)$ coincide with the internal indices of the corresponding edges of $(G,\mu)$.
In that case the scale assignment of the resulting graph is obtained in the following way. The scale assignment for the internal edges of $g$ are given by $\nu$; the scale assignment for the external edges of $g$, identified through this operation to internal edges of $G$ are given by their common value 
in $(G,\mu)$ and $(g,\nu)$.
\end{definition}

\subsection{The assigned graph combinatorial Hopf algebra}

In this section we define a  Hopf algebra on assigned Feynman graphs and we then show the relation between this structure and the combinatorics of multi-scale renormalization.

Consider now the unital associative algebra  $\cH$ freely generated by the assigned graphs, including the empty assigned graph, which we denote by $1_\cH$.

The {product} 
$m((g_1,\mu_1)(g_2,\mu_2)) = (g,\mu)$
is given by the operation of disjoint union of assigned graphs. This means that the resulting $1$PI Feynman graph  $g$ is given by the disjoint union of graphs and each disjoint component $g_i$  keeps its scale assignment $\mu_i$ ($i=1,2$) - this gives the resulting scale assignment $\mu$. 
As in the case of the Connes-Kreimer product, this product is bilinear and commutative.

 As we have already mentioned in section \ref{sec:ms}, the integers of the scale assignment $\mu$ 
are bounded by some integer cutoff $\rho$. One has $\cH_\rho\subset\cH_{\rho+1}\subset \ldots \subset
\cH_\infty$. Since we do not deal here with distinct cutoffs, we denote $\cH_\rho$ by $\cH$ in the rest 
of the paper.


\medskip

Let us now define the coproduct $\Delta:\cH\to \cH\otimes \cH$ as 
\begin{equation}
\label{coprodus} 
\Delta (G,\mu)= (G,\mu)\otimes 1_\cH+1_\cH\otimes (G,\mu)+
\sum_{(g,\nu)\subset\underline{(G,\mu)}}(g,\nu)\otimes (G/g,\mu/\nu).
\end{equation}
Note that in this definition, the high subgraphs $g$ are not necessarily connected. The coproduct can be writen explicitely as
\begin{equation}
{\Delta}(\G,\mu)=(G,\mu)\otimes 1+\sum_{(\g_{i},\nu_{i})\subset (\G,\mu)\atop \g_{i}\cap \g_{i}=\emptyset}\bigg(\prod_{i}(\g_{i},\nu_{i})\otimes(\G,\mu)\bigg)\bigg/{\textstyle \prod_{i}}(\g_{i},\nu_{i})+1\otimes(G,\mu)
\end{equation} 
where the sum runs over divergent and disjoint high 1 PI subgraphs, excluding $G$ itself.

In order to illustrate the definition of the let us list all the connected 1PI 
high superficially divergent subgraphs 
({\it i. e.} high $2-$ or $4-$point subgraphs) 
of the graph of figure \ref{scales1}
\begin{equation}
\left\{1,2,3,4\right\},\left\{7,8,9,10\right\},\left\{3,4\right\},\left\{7,8\right\}\
\end{equation}

\begin{figure}[h]
\centerline{   \includegraphics[width=4cm]{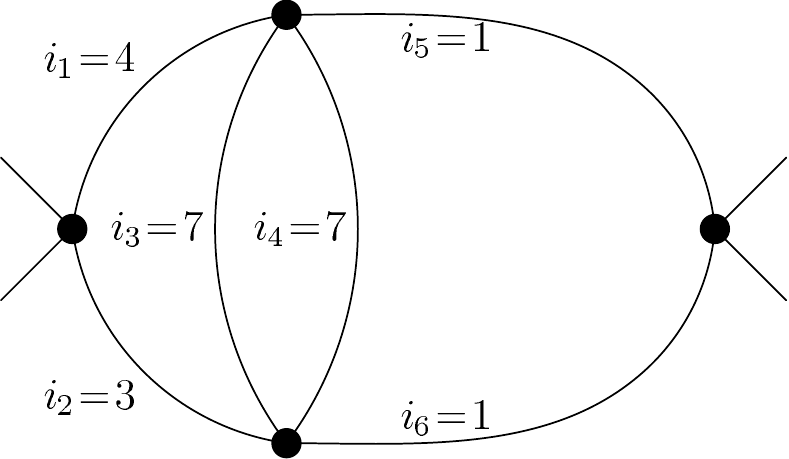}   }
\caption{The reduced graph $G\big/\left\{7,8,9,10\right\}$ }
\label{scales2}
\end{figure}

Therefore, the coproduct reads, omitting the explicit expression of the scale assignment,
\begin{multline}
\Delta(G)=G\otimes 1_\cH+1_\cH\otimes G+\left\{1,2,3,4,7,8,9,10\right\}\otimes G\big/\left\{1,2,3,4,7,8,9,10\right\}\\
+\left\{1,2,3,4\right\}\otimes G\big/\left\{1,2,3,4\right\}+
\left\{7,8,9,10\right\}\otimes G\big/\left\{7,8,9,10\right\}\\
+\left\{1,2,3,4,7,8\right\}\otimes G\big/\left\{1,2,3,4,7,8\right\}+
\left\{3,4,7,8,9,10\right\}\otimes G\big/\left\{3,4,7,8,9,10\right\}\\
+\left\{3,4\right\}\otimes G\big/\left\{3,4\right\}
+\left\{7,8\right\}\otimes G\big/\left\{7,8\right\}
+\left\{3,4,7,8\right\}\otimes G\big/\left\{3,4,7,8\right\}
\end{multline}
For example, the reduced graph $G\big/\left\{7,8,9,10\right\}$ is illustrated in Fig. \ref{scales2}.

Note that the vector space $\cH$ is graded, as in the usual Connes-Kreimer case, by the number of independent loops, number of edges or by the number of edges minus one.

\medskip

Let us recall the following result, holding for the combinatorial Connes-Kreimer Hopf algebra of Feynman graphs:

\begin{lemma} (Lemma $3.2$ of \cite{io-fabien})
\label{lem:coassociative1}
 Let $\G$ a $1PI$ Feynman graph. Provided
 \begin{enumerate}
 \item $\forall\G\in\Gsub,\,\forall\G'\in\underline{\G}$, 
the graph $\G/\G'$ is superficially divergent,
\label{item:condition1}
 \item $\forall\G_{1},\G_2$ such that $\G_1$ and $\G_2$ are superficially divergent, there exists gluing data such that $(\G_{1}\circ\G_{2})$ is superficially divergent,\label{item:condition2}
\end{enumerate}
  the  coproduct is coassociative
\begin{subequations}
  \begin{align}
\Delta\G=&\G\otimes 1
+1\otimes\G+\Delta'\G,\label{eq:Deltaprime}\\
\Delta'\G=&\sum_{\g\in\Gsub}\g\otimes\G/\g,
\end{align}
\end{subequations}
where we have denoted by $1$ the empty graph (the unit of the vector space freely generated by $1PI$ $\Phi^4$ Feynman graphs. Moreover, the notation $\Gsub$ stands for the set of superficially divergent subgraphs of $G$ ({\it i. e.} the two- and four-point subgraphs of $G$, which are not necessarily connected but with connected components which are 1PI).
\end{lemma}

This result (naturally) generalizes for assigned graphs:

\begin{lemma} 
\label{lema-noi}
 Let $(G,\mu)$ an assigned graph. Provided
 \begin{enumerate}
 \item $\forall (g,\nu)\in\underline{(G,\mu)},\,\forall (g',\nu')\in\underline{(g,\nu)}$, 
the assigned cograph $(\G/\G', \nu/\nu')$ is a two- or four-point high assigned graph,
 \item $\forall (g_{1},\nu_1), (g_2, \nu_2)$ such that $(g_1,\nu_1)$ and $(g_2,\nu_2)$ are two- or four-point high assigned graphs, there exists gluing data such that $(g_{1}\circ g_{2}, \nu_1\circ\nu_2)$ is 
a two- or four-point high assigned graph
\end{enumerate}
 the coproduct given by formula \eqref{coprodus} is coassociative
\end{lemma}
{\it Proof.} The proof of Lemma $3.2$ of \cite{io-fabien} generalizes in a straightforward manner. (QED)

\medskip

Furthermore, we define the counit $\varepsilon:\cH\to\KK$ as:
\begin{equation}
\varepsilon (1_{\cH}) =1,\ \varepsilon ((G, \mu)=0,\ \forall (G,\mu)\ne1_{\cH}.
\end{equation}
Finally, the antipode is given recursively by
\begin{align}
  S:\cH\to&\cH
\label{antipod}\\
(G,\mu)\mapsto&-(G,\mu)-
\sum_{(g,\nu)\in\underline{(G,\mu)}}S((g,\nu))(G/g,\mu/nu).\nonumber
\end{align}

This antipode can be computed, as the inverse for the convolution of the identity map. In the case of 
the Connes-Kreimer Hopf algebra of trees, this was done in \cite{fgb}. For the 
graph algebra defined here, an analogous 
computation leads to
\begin{equation}
\label{antipod-nr}
S(G,\mu)=\sum_{n=1}^L\sum_i c_{n-1}^{(i)} (-(G_1^{(i)},\mu_1^{(i)}))\ldots (-(G_n^{(i)},\mu_n^{(i)})),
\end{equation}
where
\begin{equation}
 \Delta'^{n}(G,\mu)=\sum_i c_n^{(i)}(G_1^{(1)},\mu_1^{(i)})\otimes\ldots\otimes 
(G_{n+1}^{(i)},\mu_{n+1}^{(i)})\label{antipode2}.
\end{equation}
and $c_n^{(i)}$ are the appropriate combinatorial factors obtained 
from the explicit coproduct computation.
Note that, from the explicit definition of the unit and of the counit map, one can prove that 
the sum over $n$ in (\ref{antipod-nr}) has a finite number of terms, equal to the 
number of independent cycles of the respective graph, denoted here by $L$. 
The non-recursive formula
 (\ref{antipod-nr}) is then equivalent to the sum over Zimmermann forests of 
high superficially divergent graphs.

\bigskip

We can now state the main result of this section:

\begin{theorem}
The quadrupole $(\cH, \Delta, \varepsilon, S)$ is a Hopf algebra.
\end{theorem}
{\it Proof.} We first prove the coassociativity of the coproduct \eqref{coprodus}, using Lemma \ref{lema-noi}. Let us first check the first condition of this lemma.  The fact that the resulting cograph has two or four external edges (the only thing to check in the usual Connes-Kreimer case) is trivial (since the shrinking does not affect the external structure of $g$, see Definition \ref{def:shrink}).

Let us now check in detail how the situation stands for the scales assignments.
We  denote by $i_{g'}$ the minimum of the scale assignments of the edges of $g'$ and $e_{g'}$ the maximum of the scale assignments of the edges of $g'$. Similarly, we denote by 
$i_{g}$ the minimum of the scale assignments of the edges of $g$ and $e_{g}$ the maximum of the scale assignments of the edges of $g$.  We also denote by 
$i_{g/g'}$ the minimum of the scale assignments of the edges of the cograph $g/g'$ and by  $e_{g/g'}$ the maximum of the scale assignments of the external edges of the cograph $g/g'$.

Since the external edges of $g'$ are internal edges of $g$, using Definition \ref{def:shrink}, 
this means that 
\beqa
i_{g/g'}(\mu)>e_{g/g'}(\mu)=e_{g}(\mu),
\eeqa
because, as mentioned above, the shrinking does not affect the external structure of $g$. We have thus checked the first condition of Lemma \ref{lema-noi}.

The second condition of Lemma \ref{lema-noi} is checked similarly, using Definition \ref{def:gluing}.
This concludes the proof of the coassociativity of the coproduct  Since $\cH$ is graded (see above), connected and from the coassociativity of the coproduct \eqref{coprodus}, the definition \eqref{antipod} of the antipode and from Lemma \ref{lema-manchon} leads to the result. (QED)

\bigskip

Just like in the Connes-Kreimer case, one has a straightforward pre-Lie algebra structure, given by the operation of insertion of assigned graphs. Antisymmetrizing this operation leads to a Lie algebra of assigned graphs.
Consider now the graded dual of the universal
enveloping algebra of this Lie structure. This gives the renormalization Hopf algebra defined in this section.

\subsection{Combinatorial Hopf algebras morphisms}
\renewcommand{\theequation}{\thesection.\arabic{equation}}
\setcounter{equation}{0}
\label{sec:morfisme}

Let us notice that the Gallavotti-Nicol\`o tree algebra  is 
isomorphic to the algebra $\cal H$.

\begin{proposition}
The algebra morphism  $\pi:\,{\cal H}\rightarrow {\cal H}_{\mathrm{GN}}$ defined on the generators by $\pi_{\mathrm{GN}}(G,\mu)=T_{(G,\mu)}$ is a Hopf algebra isomorphism.
\end{proposition}
{\it Proof.}  The proof is done by a direct verification.

\medskip

On the other hand, the Hopf algebra of Gallavotti-Nicol\`o trees is  a refinement of the Hopf algebra of rooted trees ${\cal H}_{\mathrm{RT}}$, as defined in \cite{CK0}. Indeed, for any 
Gallavotti-Nicol\`o tree $T$ let us define $\widetilde{T}$ as the rooted tree obtained by contracting all the arrows joining nodes decorated with the same graphs and removing all the decorations.

\begin{proposition} 
The algebra morphism defined on the generators of ${\cal H}_{\mathrm{GN}}$ by $\pi_{\mathrm{RT}}(T)=\widetilde{T}$ extends to a surjective Hopf algebra morphism from ${\cal H}_{\mathrm{GN}}$ to ${\cal H}_{\mathrm{RT}}$.
\end{proposition}
{\it Proof.}  The proof is done by a direct verification.  

\medskip

In \cite{CK1}, a graph renormalization Hopf algebra ${\cal H}_{\mathrm{ CK}}$ was introduced.  
The relation between this Hopf algebra and the one presented here is given by:



\begin{proposition} 
For every $\rho\in{\Bbb Z}_{+}$, the algebra morphism defined by on the generators of ${\cal H}_{\mathrm{CK}}$ by 
\begin{equation}
\pi^{\rho}_{\mathrm{CK}}(G)=\sum_{|\mu|\leq\rho}(G,\mu)
\end{equation}
extends to a Hopf algebra morphism from ${\cal H}_{\mathrm{CK}}$ to $\widetilde{{\cal H}}$ where $\widetilde{{\cal H}}$ is identical to ${\cal H}$ as an algebra but equipped with a coproduct that extracts all assigned graphs with 2 or 4 external edges, not only high subgraphs.
\end{proposition}
{\it Proof.} The proof is done by a direct calculation.

\medskip

For example, for the "sunset" graph below, the morphism formula above leads to:
\begin{multline}
\pi_{\mathrm{CK}}\bigg(\parbox{3cm}{\includegraphics[width=3cm]{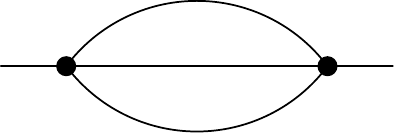}}\bigg)=
6\sum_{0\leq i_{1}<i_{2}<i_{3}\leq\rho}\parbox{3cm}{\includegraphics[width=3cm]{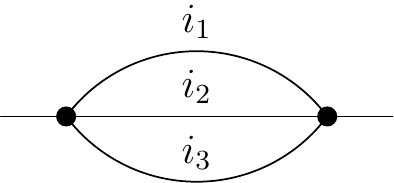}}\\
+3\sum_{0\leq i_{1}<i_{2}\leq\rho}\parbox{3cm}{\includegraphics[width=3cm]{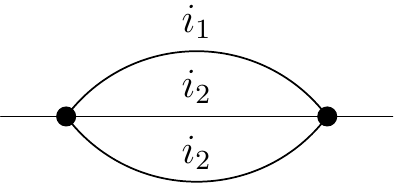}}
+3\sum_{0\leq i_{1}<i_{2}\leq\rho}\parbox{3cm}{\includegraphics[width=3cm]{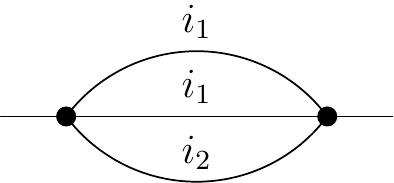}}
+\sum_{0\leq i_{1}\leq\rho}\parbox{3cm}{\includegraphics[width=3cm]{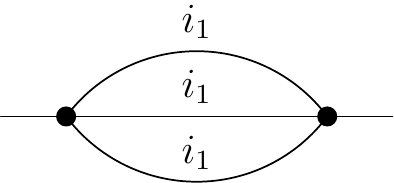}}
\end{multline}

\begin{remark}
This morphism allows us to evaluate Feynman amplitudes  with cut-off $\rho$ as
\begin{equation}
A^{\rho}(G)=A\circ\pi_{\mathrm{CK}}^{\rho}(G).
 \end{equation}
\end{remark}

\section{Multiscale renormalization combinatorics}

\subsection{Multiscale forest formula as a Hopf coaction}

\label{multiscaleHopf}
Let us now exhibit
the relation between the Hopf algebra of the previous subsection and 
the combinatorics of multiscale renormalization. 
This extends to multiscale renormalization the relation between the Connes-Kreimer Hopf algebra and renormalization.

We first recall the Feynman rules.  Given an assigned graph $(G,\mu)$, we associate a space-time variable in ${\Bbb R}^{4}$ to each vertex and and a covariance $C^{i}(x_{s(l)},y_{t(l)})$ to an edge with scale $i$ joining the vertices $s(l)$ and $t(l)$. Then, we integrate over all space-time variables but the ones attached to external edges to define the unrenormalized amplitude,
\begin{equation}
 A(G,\mu)=\int\prod_{\mathrm{ internal}\,\mathrm{ vertices}}dx_{v}\prod_{\mathrm{internal}\,\mathrm{edges}}C^{i}(x_{s(l)},y_{t(l)}).
\end{equation}
These are the Feynman rules formulated in position space, as is usual in multiscale analysis. Obviously, the evaluation of a disconnected graph is the product of the evaluation of its connected components, so that the evaluation map $ A:\,{\cal H}\rightarrow{\cal A}$ which sends an assigned graph to its value $A(G,\mu)$, with ${\cal A}$ a suitable commutative algebra depending on the variables attached $x_{1},\dots,x_{n}$ to the external edges, is a character. The connected $n$-point correlation functions  with cut-off $\rho$ are computed as a sum over over all connected Feynman graphs with $n$ external edges and scales less than $\rho$,
\begin{equation}
W(x_{1},\dots,x_{n})=\sum_{G\,\,\mathrm{connected}\,\text{assigned graph}\atop \mathrm{with}\,n\,\mathrm{external}\,\mathrm{lines}}
\sum_{|\mu|\leq\rho}A(G,\mu)[x_{1},\dots,x_{n}]\,\,
\frac{(-\lambda)^{v(G)}}{\sigma(G,\mu)}\label{connected}
\end{equation}
with $\lambda$ the coupling constant and $\sigma(G,\mu)$ the symmetry factor of the assigned graph $(G,\mu)$ (cardinal of the automorphism group of $G$ preserving the scale assignment).

In order to simplify the analytic discussion and focus on combinatorics, from now on we restrict ourselves to 
the class of $\phi^4$ graphs with at least four external edges  that do not contain any non trivial subgraph with two external edges. 
Following the terminology of \cite{book-rivasseau}, we call these graphs "biped-free graphs". Then, the multiscale 
renormalization of biped free graphs can be formulated in terms of Hopf algebras as follows.

We define  the (useful) counterterms $ C_{\mathrm{U}}$ recursively by 
\begin{equation}
 C_{\mathrm{U}}(G,\mu)=-\tau A(G,\mu)+\sum_{(G',\mu')\in\underline{(G,\mu)}}  -\tau A(G/G',\mu/\mu')\, C_{\mathrm{U}}(G',\mu'),\label{recdefC}
\end{equation}
where the sum runs over all, not necessarily connected high subgraphs whose connected components are biped-free quadrupeds (where a ``quadruped" means a graph with exactly four external edges). For any function of $n$ variables $F(x_{1},\dots,x_{n})$, $\tau F$ is defined as
\begin{equation}
\tau F(x_{1},\dots,x_{n})=\int_{({\Bbb R}^{4})^{n-1}}dx_{1}\dots dx_{n-1} F(x_{1},\dots,x_{n})
\end{equation} 
Note that if $F$ is invariant under translations, $ F(x_{1}+a,\dots,x_{n}+a)= F(x_{1},\dots,x_{n})$, $\tau F$ does neither depend on the subset of $n-1$ variables over which we integrate, nor on  the remaining variable $x_{n}$.

It is crucial to note that $ C_{\mathrm{U}}(G,\mu)$ does not depend on the space-time variable $x$ attached to the external edges because of translation invariance, so that it is a constant function.  It is also easy to see that it is multiplicative over disjoint unions and thus defines a character of ${\cal H}$. For graph with bipeds, $\tau$ involves a Taylor expansion at order 2 so that the construction is more involved.

Let ${\cal K}_{n}$ be the vector space spanned by connected assigned Feynman graphs with $n$ labeled external edges and denote by ${\cal K}$ the direct sum ${\cal K}=\oplus {\cal K}_{n}$. We define a linear map $\underline{\Delta}:\,{\cal K}\rightarrow{\cal H}\otimes{\cal K}$ by
\begin{equation}
\underline{\Delta}(\G,\mu)=(\overline{\G},\mu)\otimes 1+\sum_{(\g_{i},\nu_{i})\neq (\G,\mu)\subset (\G,\mu)\atop \g_{i}\cap \g_{j}=\emptyset}\bigg(\prod_{i}(\g_{i},\nu_{i})\otimes(\G,\mu)\bigg)\bigg/{\textstyle \prod_{i}}(\g_{i},\nu_{i})
\end{equation} 
where the sum runs over divergent and disjoint high subgraphs and $\overline{\G}$ is the Feynman graph obtained from $\G$ by erasing the labels on the external edges. $\underline{\Delta}$ is a Hopf coaction of ${\cal H}$ on ${\cal K}$, i.e. $m\circ(\epsilon\otimes\text{id})\circ\delta=\text{id}$ and 
$(\text{id}\otimes\underline{\Delta})\circ\delta=(\Delta\otimes\text{id})\circ\underline{\Delta}$. 

Furthermore, let ${\cal D}_{n}$ be a vector space of suitable distributions on  $({\Bbb R}^{D})^{n}$ in which the Feynman amplitudes take their values and ${\cal V }_{n}$ the space of linear maps from ${\cal K}_{n}$ to ${\cal D}_{n}$. The previous coaction allows us to define an action of the group of characters of ${\cal H}$ on  ${\cal V}=\oplus_{n}{\cal V}_{n}$ by $\alpha\cdot f=m\circ(\alpha \otimes f)\circ \underline{\Delta}$.

The multiscale renormalization of  biped-free graphs (see above) can then be formulated in terms of Hopf algebras by the following proposition.

\begin{proposition}
The usefully renormalized biped-free amplitudes are obtained as
\begin{equation}
A_{\text{UR}}=C_{\text{U}}\cdot A\label{multiscaleren}
\end{equation}
where the useful couterterms are defined by 
\begin{equation}
C_{\text{U}}=(\tau A)^{-1\ast}=(\tau A)\circ S
\end{equation}
\end{proposition}

This reformulation of multiscale renormalization simply relies on the fact that
\begin{equation}
S(\G,\mu)=\sum_{{\cal F}}
\prod_{(\G_{i},\mu_{i})\in{\cal F}} -(\G_{i},\mu_{i})
\end{equation}
where the sum runs over all dangerous forests, i.e. forests made of high subgraphs with four external edges.

\begin{remark}
It is crucial that $\tau A$ is a character, which follows immediately from the fact that it is a constant function. For graph with bipeds, $\tau$ involves a Taylor expansion at order 2 so that the construction is more involved and requires the introduction of new vertices of degree 2.
\end{remark}

\begin{remark}
Let us emphasize that this is not the usual BPHZ forest formula. 
Indeed, the latter involves a sum over all forests whereas the forests ${\cal F}$ 
in \eqref{multiscaleren} are such that if $G\subset G'$ in ${\cal F}$, then all edges of $G$ 
have higher scales than those of $G$ (the high scale condition, see above). 
The extra forests appearing in the BPHZ formula lead to new divergencies, 
called "renormalons" that no longer affects individual Feynman amplitudes 
but the convergence or Borel summability of the power series in the coupling constant as a whole. 
This "renormalon" problem, is cured by the multiscale expansion at the expense of 
using multiple coupling constants known as the effective coupling constants, 
as discussed in the next subsection (see \cite{book-rivasseau} for more details on this).
\end{remark}


\subsection{Effective expansion}
\renewcommand{\theequation}{\thesection.\arabic{equation}}
\setcounter{equation}{0}
\label{sec:effective}

The unrenormalized $n$-point connected correlation functions are expressed as a sum over connected Feynman graphs with $n$ labeled external edges. Thus, the bare correlation functions read 
\begin{equation}
A_{\text{bf}}(x_{1},\dots,x_{n})=\sum_{(\G,\mu),|\mu|\leq\rho\atop n\text{labeled external edges}} \frac{A(\G,\mu)[x_{1}\,\dots,x_{n}]}{\sigma(\G,\mu)}(\lambda_{\rho})^{v(\G)},\label{barecorrelations}
\end{equation}
with $v(\G)$ the number of vertices of $\G$ and $\sigma(\G,\mu)$ its symmetry factor. Note that we sum over assigned graph whose scales are bounded by $\rho$. The latter plays the role of an ultraviolet cut-off and we are ultimately interested in the limit $\rho\rightarrow\infty$. Because of the divergence of the sum over scales in \eqref{barecorrelations} when $\rho\rightarrow\infty$, one has to renormalize the Feynman graph amplitudes and expand the correlation functions into powers of the renormalized coupling constant, conventionally denoted $\lambda_{-1}$.  Then, the correlation functions read
\begin{equation}
A_{\text{bf}}(x_{1},\dots,x_{n})=\sum_{(\G,\mu),|\mu|\leq\rho\atop n\text{ labeled external edges}} \frac{A_{\text{R}}(\G,\mu)[x_{1}\,\dots,x_{n}]}{\sigma(\G,\mu)}(\lambda_{-1})^{v(\G)}.\label{renormalizedcorrelations}
\end{equation}
Here, $A_{\text{R}}(\G,\mu)[x_{1}\,\dots,x_{n}]$ denotes the renormalized Feynman graph amplitude, involving a sum over all forests, not only those made of high subgraphs. The renormalized coupling constant $\lambda_{-1}$ is computed as a sum over all  graphs with four external legs,
\begin{equation}
\lambda_{-1}(\lambda_{\rho})=
\lambda_{\rho}+\hskip-1cm\sum_{(\G,\mu),\,|\mu|\leq\rho\rho,\,i_{\G}(\mu)>i\atop\text{biped free with four external edges}}\hskip-.1cm
\frac{N(\G,\mu)}{\sigma(\G,\mu)}\,\tau A(\G,\mu)\,(\lambda_{\rho})^{v(\G)} .
\end{equation}
 As usual, $\sigma(\G,\mu)$ is the symmetry factor (cardinal of the automorphism group) while $N(\G,\mu)$ is the number of inequivalent labelings of the external edges. These numbers do not depend on the scale assignment since they involve transformations that preserve the latter.

However as discussed above, the renormalization group formalism requires to expand the correlation functions 
not in a single coupling constant $\lambda_{-1}$, but in a series of $\rho+2$ effective coupling constants 
$\lambda_{\rho},\lambda_{\rho-1},\dots,\lambda_{-1}$, one for each slice. This sequence interpolates 
between the bare coupling $\lambda_{\rho}$ and the renormalized one $\lambda_{-1}$. This is formulated in the context of Hopf algebras as follows.

For every character $\alpha$ of ${\cal H}$, let define $\rho+2$ formal power series in $\rho+2$ variables as follows
\begin{equation}
\lambda_{i}^{'}(\lambda_{\rho},\dots,\lambda_{-1})=
\lambda_{i}+\hskip-1cm\sum_{(\G,\mu),\,|\mu|\leq\rho,\,i_{\G}(\mu)>i\atop\text{biped free with four external edges}}\hskip-1cm
\frac{N(\G,\mu)}{\sigma(\G,\mu)}\,\alpha(\G,\mu)\prod_{v\,\text{vertex}}\lambda_{e_{v}(\mu)},
\quad i\in\left\{-1,\dots,\rho\right\},
\label{diffeffective}
\end{equation}
where we recall that $i_{\G}(\mu)$ is the lowest scale of the internal edges of $(\G,\mu)$ and $e_{v}(\mu)$ the highest scale on the edges attached to $v$ in $(\G,\mu)$. The inclusion of this combinatorial factor is necessary because the graphs in \eqref{diffeffective} do not carry labels on their external edges. In particular we always have $\lambda_{\rho}'=\lambda_{\rho}$ since there are no assigned graphs with $i_{\G}(\mu)>\rho$ while $\lambda'^{-1}$ involves a sum over all assigned graphs.

\begin{theorem}
The map 
 associating  the formal power series $\lambda_{i}^{'}(\lambda_{\rho},\dots,\lambda_{-1})$ to the character $\alpha$ is a group antimorphism from the group of characters $G$ of ${\cal H}$ to the group  of invertible formal power series in $\rho+2$ variables
\begin{equation}
\Psi(\beta)\circ\Psi(\alpha)=\Psi(\alpha\ast\beta). \label{morphism}
 \end{equation}
\end{theorem}

\noindent
{\it Proof:}
Let us first notice that it is sufficient to prove the result at the Lie algebra level. Indeed, any character of ${\cal H}$ can be written in a unique way as the convolution exponential of an infinitesimal character in the Lie algebra of the group of characters (see section \ref{sec:alg}). Therefore, there are infinitesimal characters $\delta$ and $\eta$ such that $\alpha=\exp_{\ast}\delta$ and $\beta=\exp_{\ast}\eta$. Then, the group morphism follows from the integration of the Lie algebra morphism using the Campbell--Baker-Hausdorff formula. Since the Lie algebra relation is linear, it is sufficient to check it for characters such that $\delta(\G,\mu)=1$ if $(\G,\mu)=(\G_{1},\mu_{1})$, $\eta(\G,\mu)=1$ if $(\G,\mu)=(\G_{2},\mu_{2})$ and vanish otherwise. 

At the infinitesimal level, the relation \eqref{morphism} reads
\begin{multline}
\sum_{\G} \frac{N(\G,\mu)}{\sigma(\G,\mu)} N((\G_{1},\mu_{1}),(\G_{2},(\mu)\mu_{2}),(\G,\mu)) \prod_{v\in V(\G)}\lambda_{e_{v}(\mu)}=\\
 \frac{N(\G_{1},\mu_{1})}{\sigma(\G_{1},\mu_{1})}\bigg(\prod_{v\in V(\G_{1})}\lambda_{e_{v}}\bigg)
 \sum_{i_{\G_{1}}(\mu_{1})>i}\frac{\partial}{\partial \lambda_{i}} 
 \bigg(
  \frac{N(\G_{2},\mu_{2})}{\sigma(\G_{2},\mu_{2})} \prod_{v\in V(\G_{2})}\lambda_{e_{v}(\mu_{2})}\bigg),
\end{multline}
where all graphs are biped free connected quadrupeds $(\G,\mu)$ such that $(\G_{1},\mu_{1})$  is a high subgraph of  $(\G_{2},\mu_{2})$ and $N((\G_{1},\mu),(\G_{2},\mu_{2}),(\G\,\mu))$ is the number of subgraphs of $(\G,\mu)$ isomorphic to $(\G_{1},\mu_{1})$ with $(\G,\mu)/(\G_{1},\mu_{1})$ isomorphic to $(\G_{2},\mu_{2})$.

Then, the results relies on the following combinatorial lemma.

\begin{lemma}
One has:
\begin{equation}
\sum_{\G} \frac{N(\G)}{\sigma(\G)} N(\G_{1},\G_{2},\G)=
 \frac{N(\G_{1})}{\sigma(\G_{1})} \frac{N(\G_{2})}{\sigma(\G_{2})} v(\G_{2})\label{combinatoriallemma}
\end{equation}

\end{lemma}
{\it Proof:}
To prove this lemma, first recall that $\frac{(4!)^{v(\G)}v(\G)!N(\G)}{\sigma(\G)}$ is the number of Wick contractions leading to the graph $\G$ in the expansion of the path integral \eqref{pathintegral}, where $N(G)$ accounts for the number of labelings of the external edges. Then,
\begin{equation}
\sum_{\G} 
\frac{(4!)^{v(\G)}v(\G)!N(\G)}{\sigma(\G)} N(\G_{1},\G_{2},\G)=
 \frac{N(\G_{1})}{\sigma(\G_{1})} \frac{N(\G_{2})}{\sigma(\G_{2})} v(\G_{2})
\end{equation}
is the number of Wick contractions leading to graphs $\G$ with a distinguished subgraph isomorphic to $\G_{1}$ such that $\G/\G_{1}$ is isomorphic to $\G_{2}$.

 Equivalently, we can start with $v(\G)=v(\G_{1})+v(\G_{2})-1$ vertices and construct $\G_{1}$. There are $\frac{(v(\G_{1})+v(\G_{2})-1)!}{v(\G_{1})!(v(\G_{2})-1)!}$ ways of choosing the vertices of $\G_{1}$ and $\frac{(4!)^{v(\G_{1})}v(\G_{1})!N(\G_{1})}{\sigma(\G_{1})}$
Wick contractions leading to $\G_{1}$. Next, we consider $\G_{1}$ as a single vertex and construct $\G_{2}$ which yields $\frac{(4!)^{v(\G_{2})-1}v(\G_{2})!N(\G_{2})}{\sigma(\G_{2})}$ Wick contractions leading to $\G_{1}$. Note that the counting involves $(4!)^{v(\G_{2}-1)}v(\G_{2})!$ instead of $(4!)^{v(\G_{2})}v(\G_{2})!$ because of the labels of the external edges of $\G_{1}$. Accordingly,
\begin{multline}
\sum_{\G} 
\frac{(4!)^{v(\G)}v(\G)!N(\G)}{\sigma(\G)} N(\G_{1},\G_{2},\G)=\\\frac{(v(\G_{1})+v(\G_{2})-1)!}{v(\G_{1})!(v(\G_{2})-1)!}
\frac{(4!)^{v(\G_{1})}v(\G_{1})!N(\G_{1})}{\sigma(\G_{1})}
\frac{(4!)^{v(\G_{2})-1}v(\G_{2})!N(\G_{2})}{\sigma(\G_{2})}
\end{multline}
which proves the lemma. (QED)

It is instructive to illustrate the combinatorics of the lemma on a simple example involving ordinary graphs. With 2 vertices, there is a single biped free quadruped, 
\begin{equation}
\parbox{2cm}{\includegraphics[width=2cm]{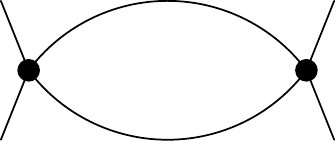}}
\qquad  \text{with}\quad\sigma(\parbox{0.8cm}{\includegraphics[width=0.8cm]{bubblelemma.pdf}})=\frac{1}{2}
\quad\text{and}\quad
N(\parbox{0.8cm}{\includegraphics[width=0.8cm]{bubblelemma.pdf}})=3
\end{equation}
$N(\parbox{0.8cm}{\includegraphics[width=0.8cm]{bubblelemma.pdf}})=3$ corresponds to the following 3 inequivalent labelings of the external edges
\begin{equation}
\parbox{2cm}{\includegraphics[width=2cm]{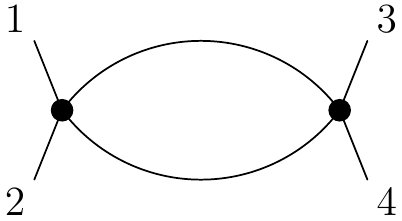}},
\qquad
\parbox{2cm}{\includegraphics[width=2cm]{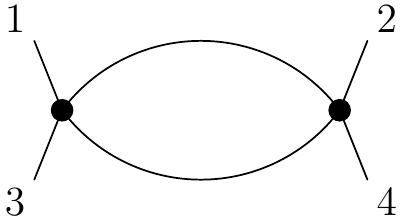}},
\qquad
\parbox{2cm}{\includegraphics[width=2cm]{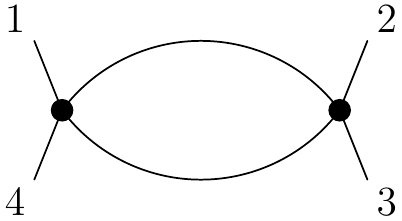}}
\qquad
\end{equation}
At order 3, we have 2 biped-free quadrupeds
\begin{equation}
\parbox{2.5cm}{\includegraphics[width=2.5cm]{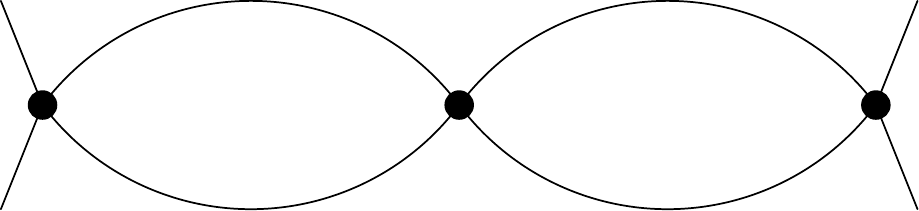}}
\;\;  \text{with}\;\; \sigma(\parbox{1.5cm}{\includegraphics[width=1.5cm]{bubblelemma5.pdf}})=\frac{1}{4},
\quad
N(\parbox{1.5cm}{\includegraphics[width=1.5cm]{bubblelemma5.pdf}})=3
\;\text{and}\;
N(\parbox{0.8cm}{\includegraphics[width=0.8cm]{bubblelemma.pdf}},
\parbox{0.8cm}{\includegraphics[width=0.8cm]{bubblelemma.pdf}},\parbox{1.5cm}{\includegraphics[width=1.5cm]{bubblelemma5.pdf}})=2
\end{equation}
and 
\begin{equation}
\parbox{1.5cm}{\includegraphics[width=1.5cm]{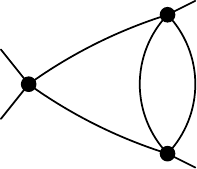}}
\qquad  \text{with}\quad\sigma(\parbox{0.8cm}{\includegraphics[width=0.8cm]{bubblelemma6.pdf}})=\frac{1}{2},
\quad 
N(\parbox{0.8cm}{\includegraphics[width=0.8cm]{bubblelemma6.pdf}})=6
\quad\text{and}\quad
N(\parbox{0.8cm}{\includegraphics[width=0.8cm]{bubblelemma.pdf}},
\parbox{0.8cm}{\includegraphics[width=0.8cm]{bubblelemma.pdf}},\parbox{0.8cm}{\includegraphics[width=0.8cm]{bubblelemma6.pdf}})=1
\end{equation}
In this case, the combinatorial lemma \eqref{combinatoriallemma} reads
\begin{equation}
\frac{3}{4}\times 2+\frac{6}{2}=\frac{3}{2}\times\frac{3}{2}\times 2 .
\end{equation}

To alleviate the notations, we have proven this lemma for ordinary graphs, not for assigned ones. In the case of assigned graphs, all goes through except that we have to take account the condition that  $(\G_{1},\mu_{1})$ is a high subgraph of $(\G_{2},\mu_{2})$, which restricts the possible insertions of $(\G_{1},\mu_{1})$ into $(\G_{2},\mu_{2})$.
(QED)

\begin{corollary}
The following power series in $\lambda_{\rho}$ are equal
\begin{equation}
\sum_{(\G,\mu),|\mu|\leq\rho\atop n\text{labeled external edges}} \hskip-0.5cm\frac{A(\G,\mu)[x_{1}\,\dots,x_{n}]}{\sigma(\G,\mu)}(\lambda_{\rho})^{v(\G)}=\hskip-0.5cm\sum_{(\G,\mu),|\mu|\leq\rho\atop n\text{ labeled external edges}} \hskip-0.5cm\frac{A_{\text{UR}}(\G,\mu)[x_{1}\,\dots,x_{n}]}{\sigma(\G,\mu)}\prod_{v\,\text{vertex}}\lambda_{e_{v}}(\mu)
\label{effectiveHopf}
\end{equation}
where the effective couplings $\lambda_{i}$ are computed using $\Psi(\tau A)$ evaluated on the bare coupling
\begin{equation}
\lambda_{i}(\lambda_{\rho})=
\lambda_{\rho}+\hskip-1cm\sum_{(\G,\mu),\,|\mu|\leq\rho,\,i_{\G}(\mu)>i\atop\text{biped free with four external edges}}\hskip-1cm
\frac{N(\G,\mu)}{\sigma(\G,\mu)}\,\tau A(\G,\mu)(\lambda_{\rho})^{v(\G)} .
\end{equation}

\end{corollary}
   
\noindent {\it Proof:} To derive this result, first compute the effective couplings $\lambda_{i}$ in terms of $\lambda_{\rho}$ using the morphism $\Psi(\tau A)$. Then, substituting the effective couplings $\lambda_{i}$ in terms of $\lambda_{\rho}$ on the RHS amounts to an action of $\tau A$. However, the usefully renormalized amplitude are precisely obtained by an action of the useful counterterms $C_{\text{U}}=(\tau A)^{-1\ast}$. Thus, the action of $(\tau A)^{-1\ast}$ due to renormalization precisely cancels the action of $\tau A$ due to the change of coupling constants.
(QED)

\begin{remark}

The counterterms defined by $C_{\text{U}}=(\tau A)^{-1\ast}$ correspond to a given renormalization scheme which amounts to Taylor subtraction at zero momentum in Fourier space.  This procedure renders the Feynman graph amplitude finite but this goal may be achieved by any other prescription. Indeed, at each step of the recursive definition of the counterterms, one can add a finite contribution $\alpha(G,\mu)$ to each Feynman graph amplitude. This amounts to transform  the counterterm as $C_{\text{U}}\rightarrow \alpha\ast C_{\text{U}}$ which in turn may be compensated by the change of effective couplings induced by $\Psi_{\alpha}$.

\end{remark}

\bigskip

\noi
{\bf Acknowledgment:}  
V. Rivasseau acknowledges Perimeter Institute grants and the ANR LQG09 grant.
 T. Krajewski and A. Tanasa acknowledge the Univ. Paris 13, Sorbonne Paris Cit\'e  
 A. Tanasa also acknowledges the grants PN 09 37 01 02 and CNCSIS Tinere Echipe 77/04.08.2010.

\noindent
{\small ${}^{a}${\it Centre de Physique Th\'eorique, 
CNRS UMR 7332, Aix Marseille Univ., Campus de Luminy, Case 907, 13288 Marseille cedex 9, France}} \\
{\small ${}^{b}${\it Laboratoire de Physique Th\'eorique,
Universit\'e Paris 11,}} \\
{\small {\it 91405 Orsay Cedex, France, EU}}\\
{\small {\it Perimeter Institute for Theoretical Physics, 31 Caroline St. N, ON, N2L 2Y5, 
Waterloo, Canada}}\\
{\small ${}^{c}${\it 
Horia Hulubei National Institute for Physics and Nuclear Engineering,\\
P.O.B. MG-6, 077125 Magurele, Romania, EU}}\\
{\small {\it Universit\'e Paris 13, Sorbonne Paris Cit\'e, \\
LIPN, Institut Galil\'ee, 
CNRS UMR 7030, F-93430, Villetaneuse, France, EU}}

\begin{thebibliography}{99}

\bibitem{flajolet}
P. Blasiak and P. Flajolet
``Combinatorial Models of Creation-Annihilation'', S\'eminaire Lotharingien de Combinatoire B65c (2011), 78 pp.

\bibitem{5}
G. Duchamp, L. Poinsot, A. Solomon, K. Penson, P. Blasiak, A. Horzela, 
``Ladder Operators and Endomorphisms in Combinatorial Physics'',  Discrete Mathematics and Theoretical Computer Science 12, 23-45, 2010.




\bibitem{CK0}
  A.~Connes and D.~Kreimer,
  {\it "Hopf algebras, renormalization and noncommutative geometry"},
  Commun.\ Math.\ Phys.\  {\bf 199} (1998) 203
  [arXiv:hep-th/9808042]


\bibitem{cm}
A. Connes and M. Marcolli, ``Noncommutative geometry, quantum fiellds and motives'', World Scientific, 2008.


\bibitem{CK1}
  A.~Connes and D.~Kreimer,
  {\it "Renormalization in quantum field theory and the Riemann-Hilbert  problem.
  I: The Hopf algebra structure of graphs and the main theorem''},
  Commun.\ Math.\ Phys.\  {\bf 210} (2000) 249
  [arXiv:hep-th/9912092]

\bibitem{CK2}
  A.~Connes and D.~Kreimer,
  {\it ``Renormalization in quantum field theory and the Riemann-Hilbert  problem.
  II: The beta-function, diffeomorphisms and the renormalization  group''},
  Commun.\ Math.\ Phys.\  {\bf 216} (2001) 215
  [arXiv:hep-th/0003188]


\bibitem{Wil1} K. Wilson, ``The Renormalization group (RG) and critical phenomena", Physical Review B, volume 4, (1971), 3174.

\bibitem{Wil2} K. Wilson, ``Renormalization group methods", Advances in Mathematics
{\bf 16}, (1975), 170-186; Nobel Lecture 1982.

\bibitem{Erice}
``Constructive quantum field theory",
Ettore Majorana International School of Mathematical Physics, eds. G. Velo and A. Wightman, 
Lecture Notes in Physics 25, Springer Verlag (1973). 

\bibitem{book-rivasseau}
V. Rivasseau, ''From perturbative to constructive renormalization``, Princeton University Press, 1991.

\bibitem{FMRS}J. Feldman, J. Magnen, V. Rivasseau and R. 
S{\'e}n{\'e}or, ``Bounds on completely convergent Euclidean Feynman Graphs", Commun. Math. Phys. {\bf 98}, 273 (1985;
and ``Bounds on renormalized Feynman graphs",  Commun. Math. Phys. {\bf 100}, 23 (1985).

\bibitem{GN}
G. Gallavotti and F. Nicol\`o, ``Renormalization theory in four-dimensional scalar fields (I) and (II),
Commun. Math. Phys. {\bf 100}, 545-590 (1985) and  {\bf 101} 247-282, (1985).

\bibitem{dCR} C. de Calan and V. Rivasseau, ``Local existence of the Borel transform in Euclidean $\phi^{4}_{4}$", 
Commun. Math. Phys. {\bf 82}, 69 (1981).


\bibitem{tHooft} G. 't Hooft, ``Rigorous construction of planar diagram field theories in four dimensional Euclidean space",
Commun. Math. Phys. {\bf 88}, (1983), 1-25.



\bibitem{Riv} V. Rivasseau, ``Construction and Borel summability of planar 4 dimensional 
Euclidean field theory", Commun. Math. Phys. {\bf 95}, 445-486 (1984).


\bibitem{BG1} G. Benfatto and G. Gallavotti
Perturbation theory of the Fermi surface
in a quantum liquid. A general quasi-particle formalism
and one dimensional systems, {\it Journ. Stat. Physics}
{\bf 59}, 541, 1990.

\bibitem{FT1} J. Feldman and E. Trubowitz,
Perturbation Theory for Many Fermions Systems, 
{\it Helv. Phys. Acta} {\bf 63}, 156 (1990).


\bibitem{Rivasseau:2011ri} 
  V.~Rivasseau,
  ``Introduction to the Renormalization Group with Applications to Non-Relativistic Quantum Electron Gases,''
  Lect.\ Notes Math.\  {\bf 2051}, 1 (2012)
  [arXiv:1102.5117 [math-ph]].
  
\bibitem{GW}
V.~Rivasseau, F.~Vignes-Tourneret and R.~Wulkenhaar,
 ``Renormalization of noncommutative phi**4-theory by multi-scale analysis,''
  Commun.\ Math.\ Phys.\  {\bf 262} (2006) 565
  [hep-th/0501036].


\bibitem{propa}
  R.~Gurau, V.~Rivasseau and F.~Vignes-Tourneret,
 ``Propagators for noncommutative field theories,''
  Annales Henri Poincare {\bf 7} (2006) 1601
  [hep-th/0512071].

\bibitem{4men}
R.~Gurau, J.~Magnen, V.~Rivasseau and F.~Vignes-Tourneret,
  ``Renormalization of non-commutative phi(4)**4 field theory in x space,''
  Commun.\ Math.\ Phys.\  {\bf 267} (2006) 515
  [hep-th/0512271].

\bibitem{gn}
F.~Vignes-Tourneret,
``Renormalization of the Orientable Non-commutative Gross-Neveu Model,''
  Annales Henri Poincare {\bf 8} (2007) 427
  [math-ph/0606069].

\bibitem{GMRT}
  R.~Gurau, J.~Magnen, V.~Rivasseau and A.~Tanasa,
``A Translation-invariant renormalizable non-commutative scalar model,''
  Commun.\ Math.\ Phys.\  {\bf 287} (2009) 275
  [arXiv:0802.0791 [math-ph]].

\bibitem{BGR} 
  J.~Ben Geloun, V.~Rivasseau and V.~Rivasseau,
 ``A Renormalizable 4-Dimensional Tensor Field Theory,''
  arXiv:1111.4997 [hep-th].

\bibitem{COR} 
  S.~Carrozza, D.~Oriti and V.~Rivasseau,
 ``Renormalization of Tensorial Group Field Theories: Abelian U(1) Models in Four Dimensions,''
  arXiv:1207.6734 [hep-th].

\bibitem{Rivasseau:2011hm} 
  V.~Rivasseau,
 ``Quantum Gravity and Renormalization: The Tensor Track,''
  AIP Conf.\ Proc.\  {\bf 1444}, 18 (2011)
  [arXiv:1112.5104 [hep-th]].
  
\bibitem{Bonzom:2011ev} 
  V.~Bonzom, R.~Gurau and V.~Rivasseau,
 ``The Ising Model on Random Lattices in Arbitrary Dimensions,''
  Phys.\ Lett.\ B {\bf 711}, 88 (2012)
  [arXiv:1108.6269 [hep-th]].
  
  \bibitem{BGS} 
V. Bonzom, R. Gurau and M. Smerlak, ``Universality in p-spin glasses",
 arXiv:1206.5539 [cond-mat.dis-nn].
  
\bibitem{Collins}
  J.~C.~Collins,
  ``Renormalization. An Introduction To Renormalization, The Renormalization Group, And The Operator Product Expansion,''
  Cambridge University Press ( 1984) 


\bibitem{io-fabien}
  A.~Tanasa and F.~Vignes-Tourneret,
 ``Hopf algebra of non-commutative field theory,''
  arXiv:0707.4143 [math-ph].

\bibitem{io-kreimer}
  A.~Tanasa and D.~Kreimer,
``Combinatorial Dyson-Schwinger equations in noncommutative field theory,''
J. Noncommutative Geometry (in press), 
  arXiv:0907.2182 [hep-th].

\bibitem{fotini}
F.~Markopoulou,
``Coarse graining in spin foam models,''
  Class.\ Quant.\ Grav.\  {\bf 20} (2003) 777
  [gr-qc/0203036].

\bibitem{io-sf}
  A.~Tanasa,
  ``Algebraic structures in quantum gravity,''
  Class.\ Quant.\ Grav.\  {\bf 27} (2010) 095008
  [arXiv:0909.5631 [gr-qc]].




\bibitem{manchon}
D. Manchon, ``Hopf algebras, from basics to applications to renormalization'', Rencontres Math\'ematiques de Glanon, 2003.

\bibitem{Kassel}
Christian Kassel, ``Quantum Groups'', 
volume 155 of Graduate Texts in Mathematics.
Springer-Verlag, 1995.

\bibitem{Dascalescu}
Sorin D\u asc\u alescu, Constantin N\u ast\u asescu, and \,Serban R\u aianu, ``Hopf Algebras,
An Introduction'', volume 235 of Pure and applied mathematics. CRC,
2001.


\bibitem{patras}
K.~Ebrahimi-Fard and F.~Patras,
  ``Exponential renormalization,'' Annales Henri Poincare 11:943-971, 2010, 
  arXiv:1003.1679 [math-ph].
``Exponential Renormalization II: Bogoliubov's R-operation and momentum subtraction schemes,''
  arXiv:1104.3415 [math-ph].

\bibitem{fgb} 
  H.~Figueroa and J.~M.~Gracia-Bondia,
 ``On the antipode of Kreimer's Hopf algebra,''
  Mod.\ Phys.\ Lett.\ A {\bf 16}, 1427 (2001)
  [hep-th/9912170].



\end{thebibliography}
\end{document}